\documentclass[smallextended,times]{article}       % onecolumn (second format)

\usepackage{graphicx}

\usepackage{mathptmx}      % use Times fonts if available on your TeX system

\usepackage{amssymb}
\usepackage{amsmath}
\usepackage{amscd}
\usepackage{amsthm}
\usepackage[all]{xy}
\usepackage{manfnt}

\usepackage{centernot}

\usepackage{bm}
\usepackage[OT2,T1]{fontenc}
\newtheorem{theo}{Theorem}[section]
\newtheorem{prop}[theo]{Proposition}
\newtheorem{lem}[theo]{Lemma}
\newtheorem{cor}[theo]{Corollary}
\theoremstyle{definition}
\newtheorem{defi}[theo]{Definition}
\theoremstyle{remark}
\newtheorem*{remark}{Remark}
\newtheorem{example}[theo]{Example}

\newcommand{\bthe}{\begin{theo}}
\newcommand{\ble}{\begin{lem}}
\newcommand{\bpr}{\begin{prop}}
\newcommand{\bco}{\begin{cor}}
\newcommand{\bde}{\begin{defi}}
\newcommand{\ethe}{\end{theo}}
\newcommand{\ele}{\end{lem}}
\newcommand{\epr}{\end{prop}}
\newcommand{\eco}{\end{cor}}
\newcommand{\ede}{\end{defi}}

\newcommand \Br {{\rm{Br}}}

\newcommand \Pic {{\rm {Pic}}}
\newcommand \val {{\rm {val}}}

\newcommand \Spec {{\rm{Spec\,}}}
\newcommand \ev {{\rm{ev}}}
\newcommand \Hom {{\rm {Hom}}}

\newcommand \rank{{\rm rank}}

\newcommand \ov{\overline}

\newcommand \Z {{\mathbb Z}}
\newcommand \Q {{\mathbb Q}}

\def\O{{\cal O}}

\newcommand\sA{{\mathcal A}}

\newcommand\beq{\begin{equation} \label}

\DeclareSymbolFont{cyrletters}{OT2}{wncyr}{m}{n}
\DeclareMathSymbol{\Sha}{\mathalpha}{cyrletters}{"58}

\usepackage{tikz-cd}

\begin{document}

\title{Brauer-Manin obstruction for zero-cycles on certain varieties }
\date{}



\author{Evis Ieronymou }


%\institute{F. Author \at
%              first address \\
%              Tel.: +123-45-678910\\
%              Fax: +123-45-678910\\
%              \email{fauthor@example.com}           %  \\


\maketitle

\begin{abstract}
We investigate the question of whether the existence of a family of local zero-cycles of degree $d$ orthogonal to the Brauer group implies the non-emptiness of the Brauer-Manin set for certain varieties.
We provide various examples of Brauer-Manin obstruction to the existence of zero-cycles of appropriate degrees.
\end{abstract}

\section{Introduction}

The theory of the Brauer-Manin obstruction for the Hasse principle is by now well developed. The parallel story about zero-cycles is not that well developed.
There are three big conjectures in this field.
The first conjecture was formulated by Colliot-Th\'el\`ene in \cite{CT03}
 and states that if $X$ is rationally
connected, then $X(k)$ is dense in $X(\mathbb A_k)^{\Br(X)}$. This generalises the same conjecture for geometrically rational surfaces by Colliot-Th\'el\`ene and Sansuc.
The second conjecture is by Skorobogatov in \cite{Skconj} and states that the Brauer-Manin obstruction to the Hasse principle is the only one for $K3$ surfaces.
The third conjecture is that the Brauer-Manin obstruction to the existence of a zero-cycle of degree $d$ is the only one \textit{for any} smooth, projective geometrically integral variety over a number field. This conjecture was put forward in various forms and for various classes of smooth, projective,
geometrically irreducible varieties by Colliot-Th\'el\`ene and Sansuc \cite{CT.S} (see also \cite{CT95}) and by Kato and Saito \cite{KS}.
We refer the
reader to \cite{W} for more information and a more refined form of the conjecture.

In this note we explore the following question. Suppose that $X/F$ is a counterexample to the Hasse principle explained by the Brauer-Manin obstruction.
What can we say about the Brauer-Manin obstruction to existence of zero cycles of degree $d$, for various $d$?
 In general not much can be said, and this is quite possibly a hopeless question. In special geometric situations however, these is something meaningful to be said. For
 example, as predicted by a result of Amer and Brumer and the conjecture on zero-cycles of degree $1$,
it is known that if $X/F$ is a del Pezzo surface of degree $4$ then a Brauer-Manin obstruction to the Hasse principle implies  a Brauer-Manin obstruction to the Hasse principle for all odd degree extensions (see Remark (1) after the Theorem).
%We study the aforementioned question for special varieties.
%To state our main result concisely we introduce the following terminology.

Let $X/F$ be a smooth, projective and geometrically integral variety over a number field, with $X(\mathbb A_F)\neq \emptyset$.
For $d\in \mathbb Z$ we  call Hypothesis (d) the following statement.

Hyp (d): There exists a family of local zero-cycles of degree d which is orthogonal to $\Br(X)$.

We also call (*) the following statement.

(*): There exists a family of local rational points which is orthogonal to $\Br(X)$

In other words, Hyp(d) is the statement that there is no Brauer-Manin obstruction to the existence of zero-cycles of degree $d$ and (*) is the statement that the Brauer-Manin set in non-empty.
In this paper we discuss when Hyp (d) implies (*). Note that if $X$ has a global zero-cycle of degree $n$ then Hyp (d) implies Hyp ($an+bd$) for any integers $a,b$.
%(cf.  proof of Prop.  \ref{bprop}).

Our main  result is the following.

\textbf{Theorem}

Let $F$ be a number field,  and $X/F$ a smooth, projective and geometrically integral variety with $X(\mathbb A_F)\neq \emptyset$.

\textit{
 \begin{itemize}
      \item[(i)] Let $d$ be an odd integer and let $X$ be a del Pezzo surface of degree 4 or a Ch\^atelet surface.  Then Hyp (d) implies (*).
 \item[(ii)]
  Let $d$ be an integer coprime to $3$,  and let $X$ be a cubic surface.  Then Hyp (d) implies (*) (\cite[Thm 1.1]{RV}).
\item[(iii)]
Let $X$ be a rationally connected variety and denote by $T$ the finite set of places consisting of the archimedean places and the finite places over which $X_{F_v}$ does not admit a smooth proper model with a separably rationally connected special fibre.
Suppose that
$X_{F_v}$ satisfies $(wP_1)$ for all $v\in T$ and  that there exists a zero-cycle of degree $n$ on $X$.  Let $d$ be an integer coprime to $n$.  Then Hyp (d) implies (*).
\end{itemize}
}

We define property $(wP_1)$ which appears in (iii) at the beginning of \textsection $3$.

\begin{remark}

\begin{itemize}
\item[(1)]
The fact that for a del Pezzo surface of degree $4$ a Brauer-Manin obstruction to the Hasse principle implies  a Brauer-Manin obstruction to the the Hasse principle for any odd degree extension was communicated to us by Colliot-Thélène.  We formalised the argument to cover more cases (see Proposition \ref{bprop}).   The key property that allows us to show that Hyp (1) implies (*) is encoded in Property $(wP_1)$.
The birational invariance of this property although is easy to show,  is not a complete triviality. There are clearly varieties beyond dimension 2 that satisfy this property, for example projective spaces.
The above might have been well-known to experts but we could not find a convenient reference.  In any case we believe that part (iii) is genuinely new and of interest.

\item[(2)]
We remind the reader that the Hasse principle holds for del Pezzo surfaces of degrees other than $2,3,4$.  It is shown in the examples that
nothing analogous to the other cases can be said for del Pezzo surfaces of degree $2$.
%\item[(5)]

%It is well-known that the Brauer-Manin obstruction to the Hasse principle is the only one for Ch\^atelet surfaces (cf.  Corollary \ref{Cha}).

\item[(3)] If Hyp (d) does not hold for $X/F$ then there exists a Brauer-Manin obstruction to the Hasse principle for $X\otimes_K L$, for any $L/F$ of degree $d$.

%\item[(7)] In particular, if as conjectured the Brauer-Manin obstruction to the existence of zero-cycles of degree $1$ is the only one,  then a del Pezzo surface of index $1$ over a number field has a rational point.

\end{itemize}
\end{remark}

In the final section we take from the literature various known counterexamples to the Hasse principle and show that
 there is a Brauer-Manin obstruction to the existence of zero-cycles of appropriate degrees.  In some cases this is not a trivial task and to the best of our knowledge such examples are quite rare in the literature.  In the case of K3 surfaces, the examples are smooth complete intersections of type $(2,3)$ and $(2,2,2)$.
It is not clear what to expect about $K3$ surfaces that are smooth complete intersections of the above kind.  and this seems to be worthy of further investigation.
  We provide some non-surface examples as well, and % We believe the techniques of the proofs would also be of interest.
one of the examples settles a question of Coray and Manoil appearing in \cite{CorayM} (see Remark after Example \ref{K3.example.CM}).

\section{Generalities  }
In this section $F$ is a number field.
We start with a small general observation.

\ble\label{obs}

\begin{enumerate}
Let $X,Y$ be smooth,  projective, geometrically integral varieties over a number field $F$.
\item[(i)] Let $f:X\to Y$ be a morphism.
If there is a Brauer-Manin obstruction to the existence of zero-cycles of degree $d$ on $Y$ then
there is a Brauer-Manin obstruction to the existence of zero-cycles of degree $d$ on $X$.

\item [(ii)] Let $f:X\dashrightarrow Y$ be a dominant rational map.
  We have an induced map $f^*:\Br( F(Y))\to \Br( F(X))$. Suppose that the image of $\Br(Y)$ under $f^*$ is contained in $\Br(X)$.

  If there is a Brauer-Manin obstruction to the existence of zero-cycles of degree $d$ on $Y$ then
there is a Brauer-Manin obstruction to the existence of zero-cycles of degree $d$ on $X$.

\item [(iii)] Let $f:X\to Y$ be a birational morphism.
 Then there is a Brauer-Manin obstruction to the existence of zero-cycles of degree $d$ on $Y$ if and only if
there is a Brauer-Manin obstruction to the existence of zero-cycles of degree $d$ on $X$.
\end{enumerate}

\ele
\begin{proof}
  Part (i) follows from functoriality.  For part (ii) let $(z_v)$ be a family of local zero-cycles of degree $d$ orthogonal to $\Br(X)$ and let $U$ be the largest open subset of $X$ where $f$ is defined.
   We can assume that for each place $v\in\Omega_F$
  the support of $z_v$ is contained in the open $U_{F_v}$. Then $f_* (z_v)$ is a family of local zero-cycles of degree $d$ on $Y$.  By our assumptions and functoriality it is orthogonal to $\Br(Y)$, which is a contradiction.
  Part (iii) follows from the other two parts.
\end{proof}

The next result concerns rationally connected varieties.
\bpr \label{dP.ob}

Let $X/F$ be a rationally connected variety and $\sA \in \Br(X)$.  Denote by $T$ the finite set of places consisting of the archimedean places and the finite places over which $X_{v}$ does not admit a smooth proper model with a separably rationally connected special fibre.
Suppose that $X(\mathbb A_F)\neq \emptyset$ and $X(\mathbb A_F)^{\sA}=\emptyset$.
Let $n$ denote the order of $\sA$ in $\Br(X)/\Br_0(X)$, and let $d\in \Z$ be coprime to $n$.
Suppose that $\sA\in f_v^{-1}(\Br_0(X_v))$ for all places $v\in T$,  where
$$
f_v:\Br(X)\to\Br(X_v)
$$
is the natural map.

Then $\sA$ gives an obstruction to the existence of zero-cycles of degree $d$ on $X$

\epr
\begin{proof}
 %Note that the smooth specialization of a del Pezzo is a del Pezzo.
We remind the reader of some notation. If $X/L$ is a variety over a field $L$, then $\Br_0(X)\subseteq \Br(X)$ is the image of the natural map $\Br(L)\to \Br(X)$ induced by the structure morphism $X\to \Spec(L)$.
If $L$ is a number field and $v$ is a place of $L$, then $L_v$ denotes the complection of $L$ with respect to $v$ and $X_v$ (or $X_{L_v}$) denotes the base extension $X\otimes_L L_v$.

We begin the proof by noting that \cite[IV Thm.3.11]{Kollar} implies that the set $T$ is indeed finite. Let $S/F_v$ be a finite extension, where $v$ is {\it any} place of $F$.

 We claim that

 $$
 \ev_{\sA,S}:X(S)\to \Q/\Z
 $$
 is constant.  If $v\in T$, then this is clear from our assumptions, and
 if $v\notin T$ then it follows from \cite[Prop. 7 ]{I}.

  As $X(\mathbb A_F)^{\sA}=\emptyset$ it follows that evaluating $\sA$ at any adelic point is a non-zero constant $a\in \Q/\Z$, with $na=0$.
 It is clear that evaluating $\sA$ at a zero-cycle of degree $d$ will give $0\neq da\in \Q/\Z$. This completes the proof.

\end{proof}

Our next two results are a slightly different presentation of the main result of \cite{RV},  and concern cubic surfaces.

\ble\label{new.V}
Let $X/L$ be a cubic surface,  with $X(L)\neq\emptyset$ and $L$ a local field.  Suppose that $3$ divides $|\Br(X)/\Br_0(X)|$. Then $CH_0(X)$ is generated by rational points.
\ele
\begin{proof}
By \cite[Lemma 2.1]{RV} $2CH_0(X)$ is generated by rational points. Therefore it sufices to show that $A_0(X)/2A_0(X)$ is generated by rational points.
By \cite[Prop. 5]{CT83} the group $A_0(X)$ embeds in $\Hom(\Br(X)/\Br_0(X), \Q/\Z )$.
Moreover by \cite{SD} the group $\Br(X)/\Br_0(X)$ is trivial or isomorphic to one of $\Z/2$, $\Z/2\times \Z/2$, $\Z/3$ or $\Z/3\times \Z/3$ as an abstract  group.
Note that the results of \cite{SD} are stated for number fields but the proofs for the abstract structure of $H^1(L,\Pic(\ov X))$, to which $\Br(X)/\Br_0(X)$ embeds,  are
geometric and hold for local fields as well; see also the first paragraph of \cite[pg. 458]{SD}.
%By our assumption and the results of \cite{SD} which describe the possible structure of $\Br(X)/\Br_0(X)$ as an abstract  group,
It follows that
  $A_0(X)$ must be a subgroup of $\Z/3\times \Z/3$.
Therefore $A_0(X)/2A_0(X)$  is trivial.  The result follows from this.

\end{proof}

\bco\label{cubic}
Let $X/F$ be a cubic surface with $X(\mathbb A_F)\neq \emptyset$. Let $d$ be an integer coprime to $3$. Then Hyp (d) implies (*).

%In particular  $X(\mathbb A_L)^{\Br(X_L)}=\emptyset$, for all extensions $L/F$ of degree coprime to 3.
\eco
\begin{proof}
  We suppose that $X(\mathbb A_{F})\neq \emptyset$,  $X(\mathbb A_{F})^{\Br(X)}=\emptyset$ and we want to show that there is an obstruction to the existence of zero-cycles of degree $d$ on $X$, for any integer $d$ coprime to 3. According to \cite[Lemma 3.4]{CT.P} the obstruction to the Hasse principle is given by a single element $\sA \in \Br(X)$, which can be taken to be of order 3 by \cite[Cor. 1]{SD} and the description
   of $\Br(X)/\Br_0(X)$ as an abstract  group in loc. cit.

    By the same reasoning as the proof of Proposition \ref{dP.ob} it suffices to check that

 $$
 \ev_{\sA,S}:X(S)\to \Q/\Z
 $$
 is constant, when  $S/F_v$ be a finite extension such that   $\sA\notin f_v^{-1}(\Br_0(X_v))$.

First we note that the map $
 \ev_{\sA,S}:X(S)\to \Q/\Z
 $ is either constant or its image has 3 elements. (see the last paragraph of \cite{CT.SD}).  Since $X(\mathbb A_{F})^{\sA}=\emptyset$ this implies that $\ev_{\sA,F_v}$ is constant.
%An element of $X(S)$ induces a zero-cycle of degree $[S:F_v]$ on $X_v$
Note also that since $\sA\notin f_v^{-1}(\Br_0(X_v))$ the assumptions of Lemma \ref{new.V} are satisfied for $X_v$. Since $\ev_{\sA,F_v}$ is constant and $CH_0(X_v)$ is generated by $F_v$-rational points by Lemma \ref{new.V} it follows that evaluating $\sA$ at a zero-cycle of $X_v$ only depends on the degree of the zero-cycle.
The result follows easily from this.

\end{proof}
\begin{remark}
\begin{itemize}

\item[(1)]

Initially we could prove this result only for some special classes of cubic surfaces.
After finishing this project we were informed about the recent preprint \cite{RV} which uses a novel geometric argument to prove the result in complete generality.
Using their geometric result  \cite[Lemma 2.1]{RV} we offer a slightly different proof.  The crux of the argument is of course still  \cite[Lemma 2.1]{RV}.

 % This result first appeared in  \cite[Thm 1.1]{RV}.  %We offer a slightly different proof of their main result.
\item[(2)]
We also note the following amusing fact about diagonal cubic surfaces over $\Q$.  If $X/\Q$ is such a surface, which is a counterexample to the Hasse principle
explained by the Brauer-Manin obstruction then it follows from \cite[Prop. 2]{CT.K.S} and the fact that at places of good reduction the evaluation map is constant that for any number field $F$, the set $X(\mathbb A_F)^{\Br(X_F)}$ is either empty or equal to the whole $X(\mathbb A_F)$.
\end{itemize}

\end{remark}

\ble
Let $X/F$ be a del Pezzo surface of degree $2$.
Suppose that $X(\mathbb A_{F})\neq \emptyset$ and $X(\mathbb A_{F})^{\mathcal A}=\emptyset$.
Suppose that for all places $v\in \Omega_F$ of bad reduction, and for all extensions $S/F_v$ of degree at most $7$ the map
$\ev_{\mathcal A,S}$ is constant. Then there is no zero-cycle of degree 1 on $X$.

\ele

\begin{proof}
  Suppose that there is a zero-cycle of degree 1 on $X$.  By \cite[Thm 4.1]{CT}. this implies that $X(L)\neq \emptyset$ for an extension $L/F$ of degree 1,3 or 7.
  However the argument in the proof of Proposition \ref{dP.ob} shows that $X(\mathbb A_L)^{\mathcal A}=\emptyset$ and so $X(L)=\emptyset$.

\end{proof}

\begin{remark}

  In principle, checking the assumptions is a finite task. Note that we may assume that $\mathcal A$ has order $2$.

\end{remark}

\section{Varieties with the $(wP_1)$ property}
In this section $L$ is a field of characteristic zero,  and $X/L$ is a smooth, proper, geometrically integral variety.
We will now introduce  property $(wP_1)$ which is a weakened variant of property $(P_1)$ which appears in \cite[pg. 302]{CT.C}.  We denote by $A_0(X)$ the subgroup of $CH_0(X)$ consisting of zero-cycles of degree 0.
Let $O\in X(L)$, and denote $f_O$ the map $f_O:X(L)\to A_0(X)$, $P\mapsto [P]-[O]$.
The surjectivity of $f_O$ does not depend on the choice of $O$.
We say that a smooth,  proper, geometrically integral variety $X/L$ satisfies $(wP_1)$ if $X(L)=\emptyset$
or the map $f_O$ is surjective, for some (any) $O\in X(L)$.

\ble
Property $(wP_1)$ is a birational property,  i.e. if $X,Y$ are  smooth proper geometrically integral varieties over $L$ that are birational over $L$ then
$X$ has $(wP_1)$ if and only if $Y$ has $(wP_1)$.

\ele
 \begin{proof}
   It is well known that $X(L)=\emptyset \iff Y(L)=\emptyset$.  We can therefore assume that $X(L)\neq\emptyset$.
Assume that there is a birational morphism $g:X\to Y$, and we can choose a point $P\in X(L)$. Let $Q=g(P)\in Y(L)$.
Consider the commutative diagram
\begin{center}
\begin{tikzcd}
X(L) \arrow[r] \arrow[d,"f_P"]
& Y(L) \arrow[d, "f_Q" ] \\
A_0(X) \arrow[r, "g_*" ]
& A_0(Y)
\end{tikzcd}
\end{center}
where $g_*$ is an isomorphism by \cite[Lem.  6.2 and Prop. 6.3]{CT.C}.

Recall that by Hironaka's results \cite[Ch.0 \textsection 5]{Hir} there is a sequence of blow-ups at smooth centers that resolves the indeterminancy locus of a rational map between smooth projective varieties.
Therefore in order to prove the lemma it suffices to consider the situation above and show two things:
\begin{itemize}
\item[(a)]  If $X$ has $(wP_1)$ then $Y$ has $(wP_1)$.
\item[(b)]
  If $Y$ has $(wP_1)$ and $X$ is the blow-up of $Y$ along a smooth center then $X$ has $(wP_1)$.
\end{itemize}
Case (a) follows immediately from the diagram,  while case (b) follows from the diagram and the fact that in this case the upper horizontal map in the diagram is surjective.

 \end{proof}

 \ble
   The following satisfy $(wP_1)$ over a field of characteristic zero.
   \begin{itemize}

   \item[(i)] Rational surfaces with a conic bundle strucure with invariant at most $5$ (see \cite[\textsection 1.2]{CT.C} for definitions of the relevant notions).
   \item[(ii)] Del Pezzo surfaces of degree $4$.
     \item[(iii)] Ch\^atelet surfaces.
     %\item[(iv)] Cubic surfaces with a line.

   \end{itemize}
 \ele
\begin{proof}
 Case (i) follows from \cite[Thm C]{CT.C}).  In the other cases, when there is a rational point they are birational to a conic bundle surface with invariant at most $5$.
 For case (ii) see \cite[Lemma 4.4]{S.Sk} and for case (iii) see \cite[Rem. 6.7 (iv)]{CT.C}.

\end{proof}

 \bpr\label{bprop}

Let $X/F$ be a variety over a number field $F$ such that $X_{F_v}$ satisfies $(wP_1)$ for all $v\in \Omega_F$ and $X(\mathbb A_F)\neq \emptyset$.
Suppose that there exists a zero-cycle of degree $n$ on $X$.  Let $d$ be an integer coprime to $n$.  Then Hyp (d) implies (*).

\epr
\begin{proof}
Let $(r_v)$ be a family of local zero-cycles of degree $d$ orthogonal to $\Br(X)$ and let $z$ denote a zero-cycle of degree $n$ on $X$.  Choose integers $a,b$ such that $ad+bn=1$.
Consider the family of local zero-cycles $(l_v)$ given by $l_v:=ar_v+bz$ for $v\in \Omega_F$.  Since $z$ is a global zero-cycle its diagonal embedding is orthogonal to $\Br(X)$ and therefore
$(l_v)$ is also orthogonal to $\Br(X)$.  Each $l_v$ is a zero-cycle of degree $1$ on $X_v$ and hence rationally equivalent to an $F_v$-point,  say $Q_v$, by the property $(wP_1)$.
The family $(Q_v)$ gives an adelic point orthogonal to $\Br(X)$ which is what we wanted to show.
\end{proof}

The following is immediate.
\bco\label{dP4CT}

Let $X/F$ be a del Pezzo surface of degree 4 with $X(\mathbb A_F)\neq \emptyset$.
Let $d$ be an odd integer.
Then Hyp (d) implies (*).

%In particular  $X(\mathbb A_L)^{\Br(X_L)}=\emptyset$, for all extensions $L/F$ of odd degree.

\eco
%\begin{proof}
  %$X$ has a zero-cycle of degree 4.
%\end{proof}

\begin{remark} According to the remarks after \cite[Lemma 3.4]{CT.P},  if there is a Brauer-Manin obstruction to the Hasse principle for a del Pezzo surface of degree 4, then the Brauer-Manin obstruction can be explained by a single element of the Brauer group of the surface of order $2$.
It is easy to see from the proof of the Proposition that the same element yields an obstruction to the existence of zero-cycles of odd degree.

\end{remark}

\bco\label{Cha}
Let $X$ be a Ch\^atelet surface with $X(\mathbb A_F)\neq \emptyset$ and $X(F)=\emptyset$.
Then there is a Brauer-Manin obstruction to the existence of zero-cycles of degree $d$ on $X$, for any odd $d$.
\eco
\begin{proof}

 A  Ch\^atelet surface has a zero-cycle of degree 4 so for $X$ we have that Hyp (d) for any odd $d$ implies (*).  Moreover, the Brauer-Manin obstruction to the Hasse principle is the only one for Ch\^atelet surfaces \cite[Thm. B]{CT.S.SD}.
The result follows from these two observations.

\end{proof}

%\begin{remark}
   %It follows trivially in this case that the existence of a zero-cycle of degree one implies the existence of a rational point .
%\end{remark}

In the case of rationally connected varieties one can use \cite[Prop. 7]{I} in order to weaken the hypothesis in Proposition  \ref{bprop}.
%for example varieties over a local field $L$ which admit a smooth model over $\mathcal{O}_L$ whose special fibre is separably rationally connected.

\bpr\label{bpropRC}

Let $X/F$ be a rationally connected variety over a number field $F$ such that $X(\mathbb A_F)\neq \emptyset$.
Let $T$ be the finite set of places containing the archimedean places and all the finite places over which $X_{F_v}$ does not admit a smooth proper model with a separably rationally connected special fibre.
Suppose that
$X_{F_v}$ satisfies $(wP_1)$ for all $v\in T$ and  that there exists a zero-cycle of degree $n$ on $X$.  Let $d$ be an integer coprime to $n$.  Then Hyp (d) implies (*).

\epr
\begin{proof}
We follow the proof of Proposition \ref{bprop} until the penultimate sentence. So we have the family $(l_v)$ which is  orthogonal to $\Br(X)$.   If $v\in T$ then we replace $l_v$ by some $Q_v$ as in the proof of loc. cit.
If $v\notin T$ then it follows from \cite[Prop. 7]{I} that we can replace $l_v$ by any $F_v$-point, call it again $Q_v$, without changing the value of the evaluation map given by any element of $\Br(X)$.
The family $(Q_v)$ gives an adelic point orthogonal to $\Br(X)$ which is what we wanted to show.
\end{proof}

\section{Examples}

\subsection{Del Pezzo surfaces}
There are many cubic surfaces with a Brauer-Manin obstruction to the Hasse principle. (see e.g. \cite[\S 7]{CT.K.S} or \cite[Ch. IV \textsection 5]{Jahnel})
\begin{example}\normalfont
 % Example from \cite{CT.K.S} and Corollary \ref{dcubQ}

  Let $X/\Q$ be the cubic surface given by
 $$
x^3+p^2y^3+pqz^3+q^2w^3=0
$$
where $p,q$ are prime numbers with $p\equiv 2 \mod 9$ and $q\equiv 5 \mod 9$

Then $X$ is a counterexample to the Hasse principle and there is a Brauer-Manin obstruction
 to the existence of zero-cycles of degree $d$ on $X$, for any integer $d$ coprime to 3.
Moreover, $X(\mathbb A_F)^{\Br(X_F)}$ is either empty or equal to the whole $X(\mathbb A_F)$, for any number field $F$.

\begin{proof}

The first assertion follows from \cite[Prop. 5, pg. 67]{CT.K.S} (and its proof) and Corollary \ref{cubic}.
The second assertion is Remark (2) after Corollary \ref{cubic}.

\end{proof}

\end{example}

The next two examples show that for del Pezzo surfaces of degree 2, anything can happen.

\begin{example}
% (Kres list)
\normalfont
Let $X/\Q$ be the del Pezzo surface of degree $2$ given by
 $$
w^2=(cz^2-x^2)(x^2+(1-c)z^2)-y^4
$$
where $c$ is a positive integer with $c \equiv 3 \mod 4$. Then there is a Brauer-Manin obstruction to the existence of zero-cycles of degree $d$, for any odd $d$.

\end{example}
\begin{proof}
  $X$ is a double cover of the Ch\^atelet surface $$Y:w^2=(cz^2-x^2)(x^2+(1-c)z^2)-y^2.$$ We have that $Y(\mathbb Q)=\emptyset$ and $Y(\mathbb A_{\mathbb Q})\neq\emptyset$ by \cite[Ex. 5.4, pg. 179]{CT.Cor.Sans}.
The result follows from Corollary \ref{Cha} and Lemma \ref{obs}.
\end{proof}

\begin{example}
% (Kres list)
\normalfont
Let $X/\Q$ be the del Pezzo surface of degree $2$ given by
 $$
w^2=-6x^4-3y^4+2z^4
$$

Then $X$ is a counterexample to the Hasse principle and we claim the following:

\begin{itemize}
  \item[(a)]
  $X(F)=\emptyset$ if $F/\Q$ has odd degree and $2$ splits completely in $F/\Q$.
  \item[(b)]
  There is no Brauer-Manin obstruction to the existence of zero-cycles of degree one.

\end{itemize}

\begin{proof}

According to \cite[Example 5]{Kres.Tchh}, we have that $X$ is a counterexample to the Hasse principle, explained by a Brauer-Manin obstruction given by
the quaternion algebra
$$
\mathcal A=(-1,h)\in \Br(X)
$$
where
$$
h=\frac{-x^2-y^2+z^2}{z^2}\in \Q(X)
$$

More precisely, $\ev_{\mathcal A}(P)$ is $0$ if $P\in X(\Q_p)$ for $p\neq 2$ and  $\ev_{\mathcal A}(P)$ is $\frac{1}{2}$ if $P\in X(\Q_2)$.
Moreover $\Br(X)/\Br_0(X)\cong \Z/2$.
% We claim that $\mathcal A$ gives an obstruction to the existence of zero-cycles of degree $d$ on $X$, for any integer $d$ coprime to 3.
To establish the claim it suffices to show that $\ev_{\mathcal A,S}$ is constant for any finite extension $S/\Q_p$, for all $p\neq 2$, and that $\ev_{\mathcal A,S}$ takes the value zero for some odd degree extension $S/\Q_2$
(cf.  the proof of Proposition \ref{dP.ob} ).

Let $P=(x_0:y_0:z_0:w_0)\in X(S)$, with $\min \{\val(x_0),\val(y_0),\val(z_0),\val(w_0)\}=0$,
where $\val$ denotes the normalised valuation of $S$. $X$ has bad reduction at $2$ and $3$, and so we need only consider the cases $p=2,3$.

{\it At the place p=2:}
Let $S=\Q_2(\pi)$, where $\pi=\sqrt[3]{2}$. An application of Hensel's lemma shows that $P=(3+7\pi^2:2+3\pi+5\pi^2:1:w_0)\in X(S)$
for some $w_0\in S$. Moreover
$$h(P)=-71\pi^2-160\pi-72\equiv \pi^2+\pi^9 \mod \pi^{10}$$
This implies that $h(P)/\pi^2 \equiv 1 \mod \pi ^7$ and hence is a square in $S$. Therefore
$$
(-1,h(P))=0
$$

{\it At the place p=3:}
We may assume that $-1$ is not a square in $S$. Looking at the equation this implies that $\val(z_0)>1$, and hence the same is true for $w_0$.
If exactly one of $x_0$ and $y_0$ has 0 valuation then clearly $\val(-x_0^2-y_0^2+z_0^2)=0$. If both of them have valuation 0 then because $-1$ is not a square we still have that
$\val(-x_0^2-y_0^2+z_0^2)=0$. In any case $\val(h(P))$ is even and this implies that $\ev_{\mathcal A,S}(P)=0$.

\end{proof}

\end{example}

\subsection{K3 surfaces and a threefold}

In this subsection the examples are $K3$ surfaces that are smooth complete intersections and a threefold that is birational to an intersection of two quadrics.

\begin{example}\label{K3.example}

\normalfont

Let $X/\Q$ be the $K3$ surface in $\mathbb P^5$ defined by
$$
X:\begin{cases}
  u^2=xy+5z^2\\
  u^2-5v^2=(x+y)(x+2y)\\
  w^2=x^2+3y^2-3z^2
\end{cases}
$$

Then $X(\mathbb A_{\Q})\neq \emptyset$ and there is a Brauer-Manin obstruction to the existence of zero-cycles of degree $d$ on $X$, for any odd $d$.

\begin{proof} This example is taken from \cite{Nguyen} where it is shown that $X(\mathbb A_{\Q})\neq \emptyset$. Moreover,
there is a morphism $X\to Y$ where $Y$ is a del Pezzo surface of degree 4, which has a Brauer-Manin obstruction to the Hasse principle (see \cite[Thm 1.2]{Nguyen} and its proof).
By Proposition \ref{dP4CT} there is a Brauer-Manin obstruction to the existence of zero-cycles of degree $d$ on $Y$, for any odd $d$.  The result now follows from Lemma \ref{obs}.
\end{proof}

\end{example}
\begin{remark}
  The same is true more generally for any $K3$ surface as in \cite[Thm 1.2]{Nguyen}. %Note that a $K3$ surface always has a zero-cycle of degree $24$.
  %Hence any of the surfaces in \cite[Thm 1.2]{Nguyen} have index $2,4$ or $8$
\end{remark}

\begin{example}\label{K3.example.CM}

\normalfont

Let $X/\Q$ be the $K3$ surface in $\mathbb P^4$ defined by
$$
X:\begin{cases}
  u^2 = 3x^2 + y^2 + 3z^2\\
5v^3 = 9x^3 + 10y^3 + 12z^3,
\end{cases}
$$

Then $X(\mathbb A_{\Q})\neq \emptyset$ and there is a Brauer-Manin obstruction to the existence of zero-cycles of degree $d$ on $X$, for any  $d$ coprime to $3$.

\begin{proof}
 This example is taken from \cite[Prop. 5.2]{CorayM} where it is shown that $X(\mathbb A_{\Q})\neq \emptyset$. There is a morphism $X\to Y$ where $Y$ is the diagonal cubic surface in $\mathbb P^3$ defined by the second equation, which has a Brauer-Manin obstruction to the Hasse principle.
By Corollary \ref{cubic} there is a Brauer-Manin obstruction to the existence of zero-cycles of degree $d$ on $Y$, for any  $d$ coprime to $3$.  The result now follows from Lemma \ref{obs}.
\end{proof}

\end{example}

\begin{remark}
 In particular $X$ has no zero-cycle of degree $1$, which settles the remark appearing right before \cite[Prop. 5.2]{CorayM}.
\end{remark}

The examples that follow do not use the results of the first threee sections.

In \cite[\textsection 6]{CT.Cor.Sans} the authors studied smooth compactifications of the variety in $\mathbb A^5$ given by
$$
u_i^2-dv_i^2=P_i(x) \quad (i=1,2)
$$
where $k$ is a field of characteristic zero, $d\in k^*$ and $d$  is not a square in $k$, and the $P_i$ are polynomials of degree $2$ in $k[x]$ with no multiple factors.
Let $Z$ be one such compactification, and assume that $k$ is a number field. They showed that $\Br(Z)/\Br_0(Z)$ is trivial  unless the $P_i$ are pairwise coprime, with all their roots in $k$.
%We remind the reader that if $\Br(Z)/\Br_0(Z)$ is trivial and $k$ in a number field then $BB_d$ holds trivially (with no content).
The next example says something in the case when $\Br(X)/\Br_0(X)\cong \Z/2$ (see \cite[Prop. 6.1]{CT.Cor.Sans}).
\begin{example}

\normalfont

Let $U/\Q$ be the  smooth threefold in $\mathbb A^5$ defined by
$$
U:\begin{cases}
  0\neq u_1^2-qv_1^2 = ax\\
0 \neq u_2^2-qv_2^2 = b(x+c)(x+d),
\end{cases}
$$
where $q$ is a prime number with $q\equiv 5 \mod8$.% (and hence $\Q_2(\sqrt{q})/\Q_2$ unramified and $-1$ is a square in $\Q_q$),

 We suppose that
\begin{itemize}
\item[(a)] $a,b,c,d\in \Q^*$,  $c\neq d$, $c,d \in \Z$,
$d$ is odd, $c$ is even,  $\val_2(a)=\val_2(b)=1$;
\item[(b)]
$\val_q(c)=1$,  $\val_q(d)\geq 2$,  $q^{-1}c$ is square in $\Q_q$,  $\val_q(a)=\val_q(b)=0$ and $a\equiv b\mod q$;

\item[(c)]
 if $p$ is a prime such that $\val_p(a)\neq 0$ or $\val_p(b)\neq 0$ or $\val_p(c-d)\neq 0$ then either $p=q$ or $q$ is a square $\mod p$.
\end{itemize}
Let $Z$ be a smooth compactification of $U$.
Suppose that  $Z(\mathbb A_{\Q})\neq \emptyset$. Then there is a Brauer-Manin obstruction to the existence of a zero-cycle of degree $d$, for any odd $d$.

\begin{proof}
 This example is taken from \cite[Prop. 7.1]{CT.Cor.Sans}, where they consider the case  $q=5$, $a=b=2$, $c=20$ and $d=25$. They show that $Z$ is a counterexample to the Hasse principle explained by the Brauer-Manin obstruction (see \cite[Rem. 7.1.4]{CT.Cor.Sans}). For the convenience of the reader we will reproduce some of the arguments in the proof of loc. cit., and we expand them a little to cover finite extensions.

 In our case we have that $\Br(Z)/\Br_0(X)$ is generated by the quaternion algebra

 $$
\mathcal A=(q,x+c)\in \Br(Z).
$$
It suffices to show that $\ev_{\mathcal A,S}$ is constant for any finite extension $S/\Q_p$, for all finite primes $p$.
We need only consider points is $U(S)$. Denote by $\val$ the normalised valuation of $S$.
Let $P\in U(S)$ and $x=x(P)$. We may assume that $q$ is not a square in $S$ and $x\neq0$.

\begin{itemize}
  \item [(i)] {\it Case $p\neq 2,q$.}
We have that $S(\sqrt{q})/S$ is an unramified extension of degree $2$. Therefore $\alpha \in S$ is a norm from $S(\sqrt{q})$
  iff $\val(\alpha)$ is even. We will show that $\val(x+c)$ is even and hence $\ev_{\mathcal A,S}(P)$ is trivial.
  By our assumptions we have that $\val(x)\equiv \val(x+c)+\val(x+d)\equiv 0 \mod 2$. If $\val(x)<0$ then $\val(x+20)=\val(x)$ and we are done.
  If $\val(x)\geq0$ then $\val(x+c)\geq 0$ and $\val(x+d)\geq 0$. Since $\val(c-d)=0$ by our assumptions either $\val(x+c)$ or $\val(x+d)$ must be zero. As their sum is even it follows in any case that $\val(x+c)$ is even.

   \item[(ii)]
{\it Case $p=2$.}
By our assumptions $S(\sqrt{q})/S$ is an unramified extension of degree $2$, and hence $\alpha \in S$ is a norm from $S(\sqrt{q})$
  if and only if $\val(\alpha)$ is even.  Let $e$ be the absolute ramification index of $S$.  We claim that
 $$e+\val(x)\equiv e+\val(x+c)+\val(x+d)\equiv 0 \mod 2.$$  Indeed,  it suffices to show that $\val(1-qt^2)$ is even for $t\in \O_S^*$.
If $\val(1-qt^2)\geq 2e+1$ then it follows from \cite[Ch. XIV,  Prop.  9]{Ser.Loc} that $q$ is a square in $S$ , contradiction. Hence we may suppose that
 $\val(1-qt^2)< 2e$.  As $1-qt^2\equiv 1-t^2\mod \pi^{2e}$, it suffices to show that $\val(1-t^2)$ is even.  If $\val(1-t)\neq\val(1+t)$ then since $e=\val(2)=\val((1+t)+(1-t))$ it would follow that
$\val(1-t^2)>2e$ which is absurd. This shows that $\val(1-t^2)$ is even and the claim is estalished.

If $e$ is odd then the argument in loc. cit.  shows that $\val(x+c)$ is odd.
 If $e$ is even we claim that $\val(x+c)$ is even.  Indeed,  otherwise we would have that $\val(x)$ is even and $\val(x+c),\val(x+d)$ odd.
 If $\val(x)\leq 0$ then it would equal $\val(x+c)$ contradiction. If $\val(x)>0$ then $\val(x+d)=0$ contradiction.

  \item[(iii)]{\it Case  $p=q$.}
Let $\pi$ be a uniformiser of $\mathcal O_S$ and write $q=\pi^eu$, for some $u\in \mathcal O_S^*$.

   Suppose that $e$ is odd.
In this case $S(\sqrt{q})/S$ is a totally tamely ramified extension of degree $2$,  and we can use \cite[Ch. V Cor. 7]{Ser.Loc}.
   Note that $q$ is a norm from $S(\sqrt{q})$ and we may choose $\pi$ so that $\pi$ is a norm from $S(\sqrt{q})$.
We claim that $x+c$ is a norm from $S(\sqrt{q})$.  First note that $a$ and $b$ are either both square or both non-squares in $S$.
If they are both non-squares then the argument is the same as in loc. cit.  If not,  the argument is similar and easier.

Suppose that $e$ is even.
 In this case $S(\sqrt{q})=S(\sqrt{u})/S$ is unramified, and  hence $\alpha \in S$ is a norm from $S(\sqrt{q})$
  if and only if $\val(\alpha)$ is even.  It is easy to check that $\val(x+c)$ is even in this case.
\end{itemize}

\end{proof}

\end{example}

Our next example is the following.
\begin{example}

\normalfont

Let $X/\Q$ be the surface in $\mathbb P^5$ defined by
$$
X:\begin{cases}
  u_1^2= xy+5v_2^2 \\
u_2^2=13x^2+950xy+32730y^2+670z^2\\
v_2^2 =-x^2-134xy-654y^2+134z^2
\end{cases}
$$
Then $X$ is a $K3$ surface with $X(\mathbb A_\Q)\neq \emptyset$,  and there is a Brauer-Manin obstruction to the existence of a zero-cycle of degree $d$, for any odd $d$.

\begin{proof}
This example is taken from \cite[Prop. 5.1]{CorayM}.  There it is shown that $X$ is smooth and has points everywhere locally.  Then they show that $X$ has no zero-cycle of degree one using \cite[Prop. 7.1]{CT.Cor.Sans} and a result of Brumer.  It is not difficult to see from our previous example that there is a Brauer-Manin obstruction to the existence of a zero-cycle of degree $d$, for any odd $d$ (as predicted by the conjecture on zero-cycles).

\end{proof}

\end{example}

\subsection{Curves}

In this subsection we confirm that for some curves that are known to have no zero-cycle of degree 1,  this absence is indeed explained by a Brauer-Manin obstruction.

\begin{example}\normalfont
  \label{bre}
Let $C/\Q$ be the curve given by
$$
C: x^4+y^4=17\cdot89z^4
$$
Then there exists a Brauer-Manin obstruction to the existence of zero-cycles of degree $d$ on $C$, for any odd $d$

\end{example}
\begin{proof}
Let $E:y^2=x^3+4\cdot17\cdot89x$ and $E':y^2=x^3-4\cdot17^2\cdot89^2x$. We have that $\rank(E)=1$ and $\rank(E')=2$.
Since $\rank(E)=1$, we see from \cite[Thm 4]{Bre} that $C$ has no points in any cubic extension of the rationals. It follows that the index of $C$ is greater than $1$ (see eg \cite[Introduction]{Cassels}). There is an isogeny from the Jacobian of $C$ to $V:=E\times E\times E'$.
  We can show using magma that $\Sha(E/\Q)[2]=\Sha(E/\Q)[4]\cong \Z/2\oplus \Z/2$ and the same is true for $E'$.
  Therefore the 2-primary part of $\Sha(V/\Q)$ is finite and hence the same is true for the  2-primary part of $\Sha(\textrm{Jac}(C)/\Q)$, see \cite[I. Lem. 7.1 (b)]{Mil} and its proof.

  We now apply \cite[Cor. 6.2.5]{Sk} (and its proof) by noting that the assumption that $\Sha(\textrm{Jac}(C)/\Q)$ is finite appearing in loc. cit.  can be weakened in our case to the assumption that the  2-primary part of $\Sha(\textrm{Jac}(C)/\Q)$ is finite. This is because the $\Q$-torsor under $\textrm{Jac}(C)$ parametrizing $0$-cycles of degree $1$ on $C$ lies in the  2-primary part of $\Sha(\textrm{Jac}(C)/\Q)$; the latter statement follows from the fact that $C$ has a rational point in an extension of degree $4$. It follows that
there is a Brauer-Manin obstruction to the Hasse principle for $C$, and that the obstruction can be given by a locally constant element $\mathcal A \in \Br(C)$. We can assume that $\mathcal A$ has order a power of $2$ since $C$ has a rational point in an extension of degree $4$. The result follows from this.

\end{proof}
\begin{remark}
More generally we have the following.
Cassels in \cite{Cassels} introduced necessary conditions for a genus $3$ curve over $\Q$ of the form
$$
C: F(x^2,y^2,z^2)=0
$$
where $F(X,Y,Z)$ is a non-singular quadratic form to have index $1$.
If these conditions are not satisfied and one can show that the Jacobians of
the genus 1 curves
$$
D_1: F(X,y^2,z^2)=0, \
D_2: F(x^2,Y,z^2)=0, \
D_3: F(x^2,y^2,Z)=0
$$ have finite 2-primary part of their Tate-Shafarevich groups, then one can conclude by the same proof as above that
there exists a Brauer-Manin obstruction to the existence of zero-cycles of degree $d$ on $C$, for any odd $d$.

\end{remark}

\begin{example}\normalfont
Let $C/\Q$ be the curve given by
$$
C: x^4+y^4=m^2z^4
$$
where $m$ is a square-free integer greate than $1$.
Suppose that the analytic rank of
$$
E:y^2=x^3+4m^2x
$$
is at most 1.

Then there exists a Brauer-Manin obstruction to the existence of zero-cycles of degree $d$ on $C$, for any odd $d$.
\end{example}
\begin{proof}
  This is similar to the proof of Example \ref{bre}, cf. the remark below it. By our assumptions the rank of $E$ is at most $1$ and so by \cite[Thm 9]{Bre} the index of $C$ is greater than $1$.
   There is an isogeny from the Jacobian of $C$ to $V:=E\times E\times E'$ where $E'$ is the elliptic curve given by $y^2=x^3-4x$. Since each factor of $V$ has analytic rank at most $1$ and since every elliptic curve over $\Q$ is modular, it follows from work of Kolyvagin that the Tate-Shafarevich group of $V$ is finite. Therefore  $\Sha(\textrm{Jac}(C)/\Q)$ is finite and we conclude like in the proof of Example \ref{bre}.

\end{proof}

Department of Mathematics and Statistics, University of Cyprus, P.O. Box 20537,
1678, Nicosia, Cyprus
\\

ieronymou.evis@ucy.ac.cy

\end{document}